\documentclass[11pt,a4paper]{article}

\usepackage[latin1]{inputenc}

\usepackage[normalem]{ulem}

\usepackage{amssymb}

\usepackage{graphicx}

\usepackage{graphicx,indentfirst,amsmath,amsfonts,amssymb,amsthm,newlfont}
\usepackage{epsfig}
\newtheorem{theorem}{Theorem}[section]

\newtheorem{definition}[theorem]{Definition}
\newtheorem{lemma}[theorem]{Lemma}
\newtheorem{corollary}[theorem]{Corollary}

\newcommand{\dif}{\mathrm{Diff}}
\newcommand{\difH}{\mathrm{Diff}(\Delta^H, M)}
\newcommand{\difV}{\mathrm{Diff}(\Delta^V, M)}
\newcommand{\Prob}{\mathbf P}
\newcommand{\R}{\mathbf R}
\newcommand{\N}{\mathbf N}

\newcommand{\F}{\mathcal F}

\newcommand{\eop }{ \hfill $\Box$ }

\begin{document}
\begin{center}


 {\Large {\bf Decomposition of stochastic flows in manifolds with complementary
distributions}}

\end{center}

\vspace{0.3cm}

\begin{center}
{\large { Pedro J. Catuogno}\footnote{E-mail:
pedrojc@ime.unicamp.br. 
Research partially supported by CNPq 302.704/2008-6, 480.271/2009-7
and FAPESP 07/06896-5.} \ \ \ \ \ \ \ \ \ \ \ \ \  Fabiano B. da
Silva\footnote{E-mail: fbsfabiano@yahoo.com.br. Research
supported by CNPQ, grant no. 142655/2005-8}

\bigskip

{ Paulo R. Ruffino}\footnote{Corresponding author, e-mail:
ruffino@ime.unicamp.br.
Research partially supported by CNPq 306.264/2009-9, 480.271/2009-7
and FAPESP 07/06896-5.}}

\vspace{0.2cm}

\textit{Departamento de Matem\'{a}tica, Universidade Estadual de Campinas, \\
13.083-859- Campinas - SP, Brazil.}

\end{center}

\begin{abstract}
Let $M$ be a differentiable manifold endowed locally with two complementary
distributions,
say horizontal and vertical. We consider the two subgroups of (local)
diffeomorphisms of
$M$ generated by vector fields in each of of these distributions.
Given a
stochastic flow
$\varphi_t$ of diffeomorphisms of $M$, in a neighbourhood of initial
condition, up to a stopping time we decompose $\varphi_t = \xi_t \circ
\psi_t$ where the first component is a diffusion in the group of horizontal
diffeomorphisms and the second component is a process in the group of vertical
diffeomorphisms.
Further decomposition will include more than two components; it leads
to a maximal cascade decomposition in local coordinates where each component
acts only
in the corresponding coordinate.
\end{abstract}


\vspace{0.3cm}
\noindent {\bf MSC2010 subject classification:} 58J65, 58D05 (57R30).

\section{Introduction}

Let $M$ be a compact differentiable manifold, we shall consider $\dif(M)$ the
Lie group of smooth diffeomorphisms of $M$ whose Lie algebra is given by smooth
vector fields. Its Lie algebra is the usual bracket operation and the
exponential map assigns to a vector field the unique flow that it generates (cf.
Hamilton \cite{Hamilton}, Milnor \cite{Milnor}) .
Given a stochastic flow $\varphi_t$
in $\dif(M)$, the decomposition of $\varphi_t$ with components in subgroups of
$\dif(M)$ which provides dynamical or geometrical information of the
system has been an interesting issue. The fact that one of the
components is again a flow (or Markovian) turns a decomposition even more
attractive in terms of applications.  In the literature this kind of
decomposition with different
aimed subgroups has appeared among others in Bismut \cite{Bismut}, Kunita
\cite{Kunita-1}, \cite{Kunita-2},  Ming Liao \cite{ML} and some of our previous
work
\cite{Ruffino}, \cite{Colonius-Ruffino}, \cite{CSR}. In the
last few papers mentioned, geometrical conditions on a Riemannian manifold have
been stated to guarantee the existence of the decomposition where the first
component lies in the subgroups of
isometries or affine
transformations. 


In this article we consider the subgroups of $\dif(M)$
whose elements preserve distributions in the sense of sections in a Grassmanian
bundle of $M$. Given a pair of distribution, say a {\it horizontal} $\Delta^H$
and a {\it vertical} $ \Delta^V$ distribution, we have associated to them the
subgroups
$\difH$ and $\difV$ generated by vector fields in the corresponding
distribution. The main question addressed here is the possibility of
decomposition of a stochastic flow as $\varphi_t = \xi_t \circ \psi_t$, with
$\xi_t \in \difH$ and $\psi_t \in \difV$. 

One of the motivation for this decomposition appears in a foliated
space, where one of the distributions is
integrable. It corresponds to one of the possible answers to the
following question: given 
trajectories of a dynamical systems in a foliated space which does not preserve
foliation, how close is the system to be leaf-preserving? Yet in another words,
how close is the vertical
component to the identity? 
Transversal perturbation in the Liouville torus in Hamiltonian systems, cf.
X.-M. Li \cite{Li}, is an example in this context. 

A principal bundle
with a connection is another natural state space in this context,
where the horizontal and vertical distributions are given by the geometry.

\bigskip

\section{Main results}


Let $M$ be a compact connected $n$-dimensional differentiable manifold, here
all geometric objects are considered smooth. Assume that (locally) $M$ is
endowed 
with a pair of regular differentiable distributions denoted
by the \textit{horizontal distribution}
 $\Delta^H:U \subset M \rightarrow Gr_k(M)$, and the \textit{vertical
distribution}
 $\Delta^{V}:U \subset M \rightarrow Gr_{n-k}(M)$, where $U\subset M$ is a
connected open set,
 $Gr_k(M)=\cup_{x\in M} Gr_k(T_xM)$ is the
Grasmannian bundle. We assume that the horizontal and the vertical
distributions are complementary in the sense that
$\Delta^H (x) \oplus \Delta^{V} (x)=T_xM$, for all
$x\in U$.

We shall consider a stochastic flow $\varphi_t$ generated by
 a Stratonovich SDE on $M$:

\begin{equation} dx_{t} \ =  \sum_{i=0}^{m}
X_{i}(x_{t})\circ dW_{t}^{i} \label{equacao original},
\end{equation}
where $W^0_t = t$, $ (W^1, \ldots , W^m)$ is a Brownian motion in $\R^m$
constructed on a filtered probability space 
$(\Omega, \F, \F_t, \Prob)$ and $X_0, X_1, \ldots , X_m$ are smooth vector 
fields in $M$. There exists a stochastic solution flow of (local)
diffeomorphisms $\varphi_t$, see e.g. among others the classical Kunita
\cite{Kunita-1}, Elworthy \cite{elworthy82}. 


In the Lie group $\dif (M)$, the
dynamics of the stochastic flow $\varphi_t$ 
is written as the right invariant equation: 

\begin{equation} d\varphi_{t} \ =  \sum_{i=0}^{m}\, R_{\varphi_t
*}
X_{i} \, \circ dW_{t}^{i} \label{equacao original no grupo de Lie},
\end{equation}
where $R_{\phi_t
*}$ is the derivative of the right translation in the group $\dif (M)$.

For a
vector field $X$ in $M$, the
associated (local) flow is denoted by $\exp \{ tX\}\in \dif(M)$. Given a
distribution $\Delta$ in $M$, we shall denote by
$\dif(\Delta, M) $ the group of diffeomorphisms
generated by exponentials
of vector fields in $\Delta$, precisely:
\begin{eqnarray*}
\dif(\Delta, M)= \mathrm{cl} \{  \exp \{t_1X_1\}\circ \ldots
\exp \{t_nX_n\}, \mbox{\  with }  X_i \in
\Delta, 
t_i \in \R,  \mbox{ for all } n\in \N \}.
\end{eqnarray*}
With this notation, the stochastic flow $\varphi_t \in \dif(TM,
M)$, a connected subgroup of $\dif(M)$ which contains the Lie subgroups $\difH$
and
$\difV$.

 If both distributions $\difH$ and $\difV$ are involutive, locally
the intersection of these two subgroups is the identity, and the elements of
each of these
subgroups preserve the 
leaves of the corresponding foliation. 


The main result of this section is the decomposition  of the stochastic
flows $\varphi_t$ in $M$ into components in
$\difH$
and $\difV$.

\begin{definition} \label{preserves_transversality} We say that a pair of
transversal complementary distributions
$\Delta^H$ and $\Delta^V$ 
preserves transversality
along the orbits of $\difH$ acting on $TM$ if
 $\xi_* \Delta^V (\xi^{-1}(x))\cap \Delta^H (x)=\{0\}$ for any element $\xi$ in
the group $\difH$. 
\end{definition}

Equivalently, $\Delta^H$ and $\Delta^V$ preserves transversality if
and only if 
 $(\xi_t)_* \Delta^V (\xi_t^{-1}(x))\cap \Delta^H (x)=\{0\}$ along trajectories
$\xi_t$ of control or stochastic systems generated by horizontal vector
fields. If $\Delta^H$ is integrable then, for any complementary distribution
$\Delta^V $, the pair $\Delta^H$ and  $\Delta^V$ preserves
 tranversality along $\difH$. In fact, $\xi_* \Delta^H = \Delta^H
\circ \xi$, hence the property follows since $\xi_*$ is an isomorphism. 

\bigskip

Our technique consists on lifting Equation (\ref{equacao original}) to
equations in the Lie subgroups 
$\difV$ and $\difH$.  

\begin{theorem}[Decomposition of flows. Global version]
\label{teorema decomposicao global}
Let $\Delta^H$ and  $\Delta^{V}$ be two complementary distributions  in
$M$ which preserves transversality along $\difH$. 
Given a stochastic flow $\varphi_t$, up to  a stopping time
there is a factorization 
\[
 \varphi_t=\xi_t \circ \psi_t
\]
where $\xi_t$ is a diffusion in $\difH$ and $\psi_t$ is a
process
in $\difV$.
\end{theorem}

\noindent \textbf{Remark 1.} The decomposition is local in the Lie group
$\difH$ and global in the manifold $M$. Compare to Remark 3.

\proof 
Denote by 
$\pi_{\Delta^V,\Delta^H}: T_xM \rightarrow \Delta^V \subset T_xM$ the projection
onto the 
vertical distribution $\Delta^V$ along $\Delta^H$.  Define, for each $i=0,1,
\ldots, m$,
and each element $\xi$ in the Lie group  $\mathrm{Diff}(\Delta^H, M)$, the
following Lie
algebra element
\begin{equation} \label{defXtilde}
 \widetilde{X_i}(x) = X_i(x)- v_i(\xi, x)
\end{equation}
where $v_i(\xi, x)$ is the unique vector in the subspace
$\mathrm{Ad}(\xi)\Delta^V\subset T_xM $ such that
$\widetilde{X_i}$ is horizontal, i.e.   
$\pi_{\Delta^V,\Delta^H} (X_i(x)-v_i(x))=0$. 


Consider the following Stratonovich SDE in the subgroup $\dif (\Delta^H, M)$
generated by 
the right action of $\xi\in \dif(\Delta^H, M)$ in the Lie algebra elements
($\xi$-dependent) $\widetilde{X}_i$:
\begin{equation} \label{equacao para a difusao}
  d\xi_t = \sum_{i=0}^{m} R_{\xi_t*} \widetilde{X}_i \  \circ dW_t^i, 
\end{equation}
with initial condition $\xi_0 = Id$. By the support theorem, the diffusion
$\xi_t$ lives
in $\difH$, since $\widetilde{X}_i$
are horizontal vector fields for all $\xi$ and all $x\in M$. Classical
existence results guarantee that there exists a solution of equation
(\ref{equacao para a difusao}) up to a stopping time. 

For the second component, we write $\psi_t= \xi_t^{-1} \circ \phi_t$ and use
that 

\[
  d \xi^{-1}_{t} = - \sum_{i=0}^{m} L_{\xi_{t}^{-1}*}\  \widetilde{X}_i
\circ
dW_t^i.
\]
 Hence, by Itô formula:

\begin{eqnarray*}
 d\psi_t &=& \sum_{i=0}^{m} \xi_t^{-1} X_i\ \xi_t \psi_t \ \circ dW_t^i -
\xi_t^{-1} \widetilde{X}_i\ \xi_t \psi_t\   \circ dW_t^i. \\
 & =&  \sum_{i=0}^{m} \mathrm{Ad}(\xi_t)(X_i -\widetilde{X}_i )\  \psi_t\  
\circ dW_t^i \\
 &=&  \sum_{i=0}^{m} \mathrm{Ad}(\xi_t) (v_i) \  \psi_t\   \circ dW_t^i.
\end{eqnarray*}
Since, by construction,  $\mathrm{Ad}(\xi_t) (v_i)  (x) \in \Delta^V (x)$, again
by the support theorem, we have 
that $\psi_t \in \mathrm{Diff}(\Delta^V, M)$. 

\eop

Although in general the factor $\xi_t$ is not a solution of an autonomous SDE in
the manifold itself, if the 
vertical distribution is invariant by $\mathrm{Ad}(\xi_t)$, for  $t \geq 0$, 
then $\xi_t$ is indeed a solution of an SDE generated by  (horizontal) vector
fields in $M$:

\begin{corollary}[Horizontal diffusion on the manifold] \label{decomposicao com
equacao na variedade}
  If $\mathrm{Ad}(\xi_t)\Delta^V = \Delta^V $ for $t \geq 0$ 
then $\xi_t$ is the solution of the following equation in $M$:
\[
  dx_{t} \ =  \sum_{i=0}^{m} X^H_{i}(x_{t})\circ dW_{t}^{i},
\]
where $X^H_ i$ are the horizontal component  $\pi_{\Delta^H, \Delta^V}X_i$, for
each $i=0, 1, \ldots, m$. 
\end{corollary}
\proof Indeed, in this case, for each $i=0, 1, \ldots, m$, the vectors $v_i$ in
the proof of Theorem \ref{teorema decomposicao global}
is simply the vertical component $X^V_ i = \pi_{\Delta^V, \Delta^H}X_i$.

\eop

\begin{corollary}[Constant energy foliation] \label{decomposicao caso
hamiltoniano}
Let $M$ be a Riemannian manifold and $h:M \rightarrow \R$ be a
submersion. Then, up to a stopping time, 
the stochastic flow $\varphi_t$ of equation (\ref{equacao original}) can be
factorized as 
\[
\phi_t=\xi_t \circ \psi_t
\]
where the diffusion $\xi_t$ preserves $h$ and $\psi_t$ is a process in the group
of
diffeomorphisms which acts orthogonally
on the leaves of constant $h$.
\end{corollary}
\proof

Take $\Delta^H = \mathrm{Ker}\  h_*$ and $\Delta^V= \{ \lambda \nabla h;
\lambda \in \R \}$. These distributions are
involutive and orthogonal.
In this particular case, the vertical vectors $v_i(x)$, $i=0,1,
\ldots, m$  in the
proof of Theorem
 \ref{teorema decomposicao global} can easily be 
calculated at each point $x\in U$ by:
\[
 v_i(x) = \frac{<\nabla h(x), X(x)>}{<\nabla h(x), \xi_* \nabla h (\xi^{-1}(x))
>} \xi_* \nabla h (\xi^{-1}(x)).
\]

\eop

\bigskip

 Given a neighbourhood $U_{x_0}$ of
 a point $x_0 \in M$, we shall denote by $\dif(\Delta^H, U_{x_0}) = \{\varphi:
U_{x_0} \rightarrow \varphi (U_{x_0}) \}$  the set of diffeomorphisms
of $U_{x_0}$ generated by horizontal vector fields; analogously, $\dif(\Delta^V,
U_{x_0})$ denotes the diffeomorphisms generated by vertical vector fields.

Transversality condition along the
orbit of $\difH$ required in Theorem \ref{teorema decomposicao global} can be
suppressed
in the local version.

\begin{theorem}[Decomposition of flows. Local version]
\label{teorema decomposicao local}
Assume that the manifold $M$ is locally endowed with a pair of complementary
distributions
$\Delta^H$ and  $\Delta^{V}$.
Given a stochastic flow of local diffeomorphisms $\varphi_t$,
for each point $x_0 \in M$ there exists an open
neighbourhood $U_{x_0}\subset M$ such that, up to a stopping time, we can
decompose 
\[
 \varphi_{t \left|_{U_{x_0}} \right. }  = \xi_t \circ \psi_t
\]
where $\xi_t$ is a diffusion in $\dif(\Delta^H, U_{x_0})$ and $\psi_t$ is a
process
in $\dif(\Delta^V, U_{x_0})$.
\end{theorem}

\proof 
The adjointly vertical correction term $v_i(\xi, x)$ in equation
(\ref{defXtilde}) depends continuously on $\xi$. At the identity
$\xi_0=Id$, this equation gives $v_i(\xi_0,
x)= \pi_{\Delta^V,\Delta^H} X(x)$, therefore it is well defined in a
neighbourhood of $\xi_0=Id$ (in the space of diffeomorphisms), and in a
neighbourhood of the initial
condition $x_0$ (in $M$), where $\Delta^H$ and  $\Delta^{V}$ are defined and
complementary. So, $v_i(\xi, \cdot)$ is defined up to a stopping time $\tau$,
such that for $t< \tau (\omega)$ at the point $\varphi_t(\omega, x_0)$, $\omega
\in \Omega$, the
distributions $\Delta^H $ and $\Delta^V$ are defined with $Ad (\xi_t) \Delta^V $
and $\Delta^V$  complementary. Hence, the equations of the horizontal and
vertical
components $\xi_t$ and $\psi_t$ respectively, hold up to the minimum of
the explosion of equation (\ref{equacao para a difusao}) and the stopping time
$\tau$.
 
\eop

\begin{corollary}[Involutive distributions] \label{decomposicao com folheacao
involutiva} If both 
$\Delta^H$ and $\Delta^V$ are involutive then the local decomposition is unique.
\end{corollary}

\proof Just note that in a neighbourhood
 $U_{x_0}\subset M$  holonomy of the foliations vanish, i.e. 
$\dif(\Delta^H,U_{x_0})  \cap
\dif(\Delta^V,U_{x_0})=
Id$.

\eop

\bigskip

\noindent \textbf{Example 1:} Let $M=\R^3$ with the canonical
basis $\{e_1, e_2, e_3\}$ and
$p=(x,y,z)\in M$. Consider the distributions
$\Delta^H (p) = \mathrm{span}
 \{ (\cos y^2, 0, \sin y^2), (0,1,0) \}$ and  $\Delta^V (p) = \mathrm{span}
 \{ (-\sin y^2, 0, \cos y^2) \}$. For the constant vector field $Y\equiv e_2$ in
$\Delta^H (x)$, the linearization
of the corresponding flow $\varphi_t \in \difH$ is the identity $d\varphi_t =
1d$. There exists a sequence of points $p_n$ and a sequence $t_n
\rightarrow 0$ such that transversality degenerates i.e. $\varphi_{t_n} \Delta^V
({\varphi_t}^{-1}(p_n)) \subset \Delta^H (p_n)$.
Hence the
pair of transversal distribution $\Delta^H$ and
 $\Delta^V$ are complementary but they do not preserve transversality
globally along $\difH$ according to Definition \ref{preserves_transversality}.
Local decomposition as in Theorem \ref{teorema decomposicao local}
holds. Note that changing the order between the horizontal and the
vertical
distributions, we have that the pair of distributions does preserve
transversality along
$\difV$ since $\Delta^V$ is integrable. 

\bigskip

\noindent \textbf{Example 2:} Coordinate systems in differentible manifolds
is a natural source of complementary involutive distributions. For
instance, with the pair of 
foliation on $\R^n-\{ 0\}$ given by spheres centred at the origin and the
corresponding 
radial lines, the angular difusion $\xi_t$ projected in $S^{n-1}$  and the
radial component $\psi_t$ of Theorem \ref{teorema decomposicao
global} is one of the examples related to Liao's factorization in
\cite{Liao-SPA} and \cite{Liao-JTP}.

\bigskip

\noindent \textbf{Example 3:} Let $(\pi:P\rightarrow M,G)$ be a principal fibre
bundle  with a connection form $\omega$. It is convenient to define $\Delta^H =
\ker
\pi_*$, the tangent subspace of the orbits of the action of $G$ and  $\Delta^V =
\ker \omega$, established by the geometry. In general a flow $\varphi_t$ in
manifold with two complementary distributions can not be decomposed into two
diffusions in $\difH$ and in $\difV$ simultaneously. 
In the case of $P$ being the frame bundle of the differentiable 
manifold $M$ with structural group $Gl(n,\R)$ we have an interesting particular
situation: For the stochastic flow $\varphi_t$ in the base manifold $M$,
consider the induced linearized flow $ \varphi_{t*}$ in $P$. It is well known
that many informations for the dynamics of $\varphi_t$ in $M$ can be studied
using the vertical
and horizontal components of $\varphi_{t*}$, e.g. parallel transport, Lyapunov
exponents, rotation number, and others (see e.g. among others
\cite{nomizu}, \cite{elworthy}, \cite{Elw-Xue-Mei}
\cite{Ruffino}). We can write (see e.g. \cite{elworthy})
 
\[
 \varphi_{t*} (u_0) = R_{g_t} \ (u^h_t),
\]
where $R_{g_t}$ is the right action of the structural group, $g_t$ obtained by
parallel transport and $u^h_t \in P$ is
the horizontal lift of an initial frame $u_0$ with $\pi(u_0)=x_0$. 
For any continuous $\gamma_t$ with $\gamma_0=Id$ in the group of holonomy of
$P$ we also have that 
\[
 \varphi_{t*} (u_0) = R_{\gamma_t g_t} \ (u^h_t \gamma_t^{-1}),
\]
which corresponds to one of the decompositions stated in Theorem \ref{teorema
decomposicao local}. If the
manifold has curvature zero then, by Corollary \ref{decomposicao com
folheacao involutiva} we have uniqueness of decomposition, hence the
horizontal component in $\dif (\Delta^H, P)$ is $\xi_t= R_{g_t}$ and $\psi_t
(u_0)=
u^h_t$. Moreover, the right action
of the structural group 
in $P$ preserves the connection $\omega$, hence  
Corollary \ref{decomposicao com equacao na variedade} says that one can obtain
an equation
for $\xi_t$ in $P$ instead of in $\dif(\Delta^H, P)$, as it is expected by the
geometry.

\bigskip

\noindent \textbf{Remark 2.}  If in
equations (\ref{equacao original no grupo de Lie}) and (\ref{equacao para a
difusao}), instead of right
translation
$R_{\varphi*}$ one considers the left translation $L_{\varphi*}$, 
one finds, left diffusions $\varphi^L_t$, $\xi^L_t$ respectively and a
vertical process $\psi_t^L$
such that the decomposition stated in the theorems above changes the order:
\[
 \varphi^L_t = \psi^L_t \circ \xi^L_t.
\]

\eop

\section{Cascade Decomposition}

In this section we assume that 
we have a sequence of
complementary distributions. 
Precisely, consider two sequences of enclosing
distributions $\Delta^H_1 \subset \ldots \subset \Delta^H_k$ and $\Delta^V_1
\supset \ldots \supset \Delta^V_k$. This is equivalent of saying that
$(\Delta^H_1,
\ldots, \Delta^H_k)$ and
$(\Delta^V_k,  \ldots \Delta^V_1)$ are 
smooth sections of a flag bundle over $M$. 
The sequence of distribution $(\Delta^H_1,
\ldots, \Delta^H_k)$ is a section of the maximal flag manifold if $\dim
\Delta^H_{i+1} - \dim
\Delta^H_{i}=1$ and $k=n$.

We assume that for
each $i=1, \ldots, k \leq n$ the pair $\Delta^H_i$ and $\Delta^V_i$ are
complementary, i.e. $\Delta^H_i(x) \oplus \Delta^V_i(x)= T_{x}M$ for every $x\in
M$. The enclosing hypothesis on the subspaces (i.e. they are
sections of a flag bundle) induces an
enclosing property in the corresponding generated Lie groups: $\dif
(\Delta^H_i,M)
\subseteq \dif(\Delta^H_{i+1},M)$ and
$\dif(\Delta^V_{i+1},M) \subseteq \dif(\Delta^V_{i},M)$ for each $i=1, \ldots ,
k-1$.

\begin{theorem}[Cascade decomposition: Global
version]
\label{decomposicao cascata global}
Let $(\Delta^H_1, \Delta^H_2,  \ldots , \Delta^H_k)$ and
$(\Delta^V_k, \Delta^V_{k-1},  \ldots , \Delta^V_1)$ be sequences of
enclosing distributions (i.e. smooth sections of flag
bundles) such
that the pairs $\Delta^H_i$ and  $\Delta^{V}_i$, $i=1, \ldots k$, are
complementary in
tangent spaces and transversality is preserved along
 the action of $\dif(\Delta^H_i, M)$. 
Given an stochastic flow $\varphi_t$ generated by equation (\ref{equacao
original no grupo de Lie}), up to a stopping time 
we can decompose 
\[
 \varphi_t=\xi^1_t \circ \ldots \circ \xi_t^k \circ \Psi_t
\]
where for each $i=1, \ldots, k$, $\xi^i_t \in \dif(\Delta^H_i, M)$,
the 
composition of the first $i$-th component $(\xi^1_t \circ \xi_t^2 \ldots \circ
\xi^i_t)$
is a diffusion in $\dif(\Delta^H_i(M)$ and the composition of the last
 components $(\xi^{i+1}_t \circ \ldots  \circ \xi^{k}_t \circ \Psi_t)$ is
a process
in $\dif(\Delta^V_{i}; M)$. 
$\Psi_t=Id$ if $\dim \Delta^H_k=n$.
\end{theorem}

\proof By Theorem \ref{teorema decomposicao global}, for each $i=1, \ldots, k$
there exists a decomposition $\varphi_t = \tilde{\xi}^{(i)}_t \circ
\tilde{\Psi}_t^{(i)}$ such that $\tilde{\xi}^{(i)}_t$ is a diffusion in
$\dif(\Delta^H_i, M)$ and $\Psi_t^{(i)}$ lives in
$\dif(\Delta^V_i, M)$. The result follows by taking $\xi^1 =
\tilde{\xi}^{(1)}_t$ and by induction 
\[
 \xi^i_t = \left( \tilde{\xi}^{(i-1)}_t \right)^{-1} \circ\ \tilde{\xi}^{(i)}_t
\]
for $1<i \leq k$ and $\Psi_t = \tilde{\Psi}_t^{(k)}$. If $\dim \Delta^H_k=n$
then $\dif(\Delta^V_k,M)= \{Id\}$, which proves the last statement.

\eop

\begin{corollary}[Cascade decomposition: local
version]
\label{decomposicao cascata local}

Assume that $M$ is locally endowed with pairs of
complementary distributions $\Delta^H_i$ and  $\Delta^{V}_i$, $i=1,2 \ldots, k$,
such that the sequences  $(\Delta^H_1, \Delta^H_2, \ldots , \Delta^H_k)$ and
$(\Delta^V_k, \Delta^V_{k-1},  \ldots , \Delta^V_1)$ are enclosed
distributions (i.e. sections of flag
bundles). 
Given the stochastic flow of local diffeomorphisms $\varphi_t$ generated by
equation
(\ref{equacao
original no grupo de Lie}), for each point $x_0 \in M$ there exists an open
neighbourhood $U_{x_0}\subset M$ such that, up to a stopping time, we can
decompose 
\[
 \varphi_{t \left|_{U_{x_0}} \right. }  = \xi^1_t \circ \ldots \circ \xi_t^k
\circ \Psi_t
\]
where for each $i=1, \ldots, k$, $\xi^i \in \dif(\Delta^H_i, U_{x_0})$, the 
composition of the first $i$-th component $(\xi^1_t \circ \xi_t^2 \circ \xi^i)$
is a diffusion in $\dif(\Delta^H_i, U_{x_0})$ and the composition of the last
components $(\xi^{i+1}_t \circ \ldots  \circ \xi^{k} \circ \Psi_t)$ is
a process
in $\dif(\Delta^V_{i}; U_{x_0})$. 
$\Psi_t=Id$ if $\dim \Delta^H_k=n$.
\end{corollary}

\proof By Theorem \ref{teorema decomposicao local}, for each $i = 1,2, \ldots,
k$ there exists a local decomposition
$\varphi_t = \tilde{\xi}^{(i)}_t \circ
\tilde{\Psi}_t^{(i)}$ which holds up to a
stopping time $\tau_i$. The result follows up to  $\tau =
\min\{\tau_1, \tau_2,
\ldots, \tau_k \}$ repeating the construction of $\xi^i_t$,  as in the proof of
the Theorem \ref{decomposicao cascata global}.

\eop

A particularly interesting situation is when the tangent space is
decomposed as a direct sum of one-dimensional subspaces and we take
maximal flag sections whose integrable distributions are generated by
direct sum of these one-dimensional subspaces.

\begin{corollary}[Decomposition preserving local coordinates]
\label{decomposicao coordenadas locais}
Let $U \subset M$ be an open set with local coordinates $\phi=(\phi_1, \ldots,
\phi_n): U \subset M \rightarrow \R^n$.
Given a stochastic flow of local diffeomorphisms $\varphi_t$ generated by
equation
(\ref{equacao
original no grupo de Lie}), up to a stopping time, we can decompose locally
\[
 \varphi_t=\xi^1_t \circ \ldots \circ \xi_t^n
\]
where for each $i=1, \ldots, n$, the (local) diffeomorphism $\xi^i_t$
preserves the $j$-th coordinates for all $j\neq i$ (i.e. $\phi_j
(\xi^i_t)(x)$ is
constant for each $x$ in the domain); the 
composition $(\xi^1_t \circ \xi_t^2 \circ \xi^i_t)$
is a diffusion of diffeomorphisms which preserves
coordinates $\phi_{i+1}, \phi_{i+2}, \ldots , \phi_n$. 


The decomposition, in this order, is unique.
\end{corollary}

\proof

Define the sequence of complementary involutive distributions by the following:
at each $x \in U$, $\Delta^H_i = \mathrm{span} \{\nabla \phi_1, \ldots, \nabla
\phi_i \} $ and $\Delta^V_i = \mathrm{span} \{\nabla \phi_{i+1}, \ldots,
\nabla \phi_n \} $, $i=1, \ldots, n$. Diffeomorphisms in $ \dif(\Delta^H_i,
M)$ preserve the leaves of the
foliation induced by $\Delta^H_i$, i.e. it preserves coordinates $\phi_j$ if
$j>i$. Analogously, diffeomorphisms in $\dif(\Delta^V_i, M)$  preserves
coordinates $\phi_j$ if
$j\leq i$.

Corollary \ref{decomposicao cascata local} guarantees that 
\[
 \varphi_t = \xi_t^1 \circ \ldots \circ \xi_t^n
\]
with $(\xi_t^1 \circ \ldots \circ \xi_t^i) \in \dif(\Delta^H_i, M)$ for $i=1,
\ldots, n$. Hence, since both $(\xi_t^1 \circ \ldots \circ \xi_t^i)$ and
$(\xi_t^1 \circ \ldots \circ \xi_t^{i-1})$  preserve coordinates
$\phi_j$ for $j>i$, then $\xi_t^i$ also preserves coordinates
$\phi_j$ for $j>i$.

Moreover, $(\xi_t^i \circ \xi_t^{i+1} \circ  \ldots \circ \xi_t^n) \in
\dif(\Delta^V_{i-1}, M)$. Hence, since both  $(\xi_t^{i+1} \circ \ldots \circ
\xi_t^n)$ and
$(\xi_t^i \circ \ldots \circ \xi_t^{n})$  preserve coordinates
$\phi_j$ for $j\leq i-1$, then $\xi_t^i$ also preserves coordinates
$\phi_j$ for $j\leq i-1$.

\eop


Next Lemma gives the picture of the restriction for existence of the
decomposition we treat in this article. This restriction appears as an
explosion time in the SDE for the components of the decomposition, cf. equation
(\ref{equacao para a difusao}).

Given a diffeomorphism $\varphi:U \rightarrow V$  between open
sets
$U, V \in \R^n$, with $x=(x_1, \ldots, x_n)\mapsto (\varphi_1(x), \ldots ,
\varphi_n(x))$. We denote the lower right
$(n-i+1)\times (n-i+1)$-submatrix of the differential $\varphi_*$ with respect
to
the canonical basis by 
\[
 \frac{\partial (\varphi_i, \ldots, \varphi_n)}{\partial (x_i, \ldots, x_n)} :=
\left[
\frac{\partial \varphi_j}{ \partial x_k} \right]_{ i\leq j\leq  n, \ i\leq k
\leq n}.
\]

Let $\mathrm{Dif}^i(U)$ be the set of diffeomorphisms $\xi: U
\rightarrow \xi(U)$ which only acts on the $i$-th coordenate, i.e.
$\xi(x)= (x_1, \ldots, x_{i-1}, \xi_i (x), x_{i+1}, \ldots, x_n)$.

\begin{lemma} \label{decomposicao total} Let $U \subset \R^n$ be an open set.

\noindent \textbf{(i)} A diffeomorphism  $\varphi:U
\rightarrow V$ can be written as a
composition of a sequence of diffeomorphisms 
\begin{equation} \label{coordenada a coordenada}
 \varphi = \xi^1 \circ \ldots \circ \xi^n 
\end{equation}
where each $\xi^i \in \mathrm{Dif}^i(U_x)$ changes only the $i$-th coordenate in
an open neighbourhood $U_x\subset U $ of $x$ if and only if
$ \det \frac{\partial (\varphi_i, \ldots, \varphi_n)}{\partial (x_i, \ldots,
x_n)}(x)
\neq 0$ for all $1\leq i \leq n$. 

\bigskip

\noindent \textbf{(ii)} The decomposition of equation (\ref{coordenada
a coordenada}), in this order, is unique.

\bigskip

\noindent \textbf{(iii)} If for a certain $1 \leq i \leq n$,  $ \det
\frac{\partial (\varphi_i, \ldots, \varphi_n)}{\partial (x_i, \ldots,
x_n)}(x)$ goes to zero as $x \rightarrow p\in \mathrm{cl}(U)$ then $||
\xi^i ||_{C^1}$ increases to infinity in $U_x$ .

\bigskip

\noindent \textbf{(iv)} The set of diffeomorphisms $\varphi \in \mathrm{Dif}(U)$
which admits the
decomposition above is a dense open set containing the identity in the space of
diffeomorphisms $\mathrm{Dif}(U)$ with respect to the  $||\cdot||_{C^1}$
topology.
\end{lemma}
\proof
(i)  Write in local coordenates $\varphi = (\varphi_1, \ldots, \varphi_n)$. The
maps
\[
x \mapsto (x_1, \ldots, x_k, \varphi_{k+1}(x),
\ldots, \varphi_n(x))
\]
 are diffeomorphisms in a neighbourhood of $x\in U$.
Define, for each $1\leq k\leq n$ the diffeomorphisms
\begin{equation}  \label{formula para decomp linear}
 \xi^k = (x_1, \ldots, x_{k-1}, \varphi_k (x), \ldots, \varphi_n(x)) \circ (x_1,
\ldots, x_k, \varphi_{k+1}(x),
\ldots, \varphi_n(x))^{-1}.
\end{equation}
They satisfy the required condition. The converse is trivial.

\bigskip


%
%

\noindent (ii) Uniqueness follows because the intersections of the sets
$\mathrm{Dif}^i(U_a)$ and $\mathrm{Dif}^j(U_a)$ is the identity if $i\neq j$.

\bigskip

\noindent (iii) It follows by formula (\ref{formula para decomp linear}).

\bigskip

\noindent (iv) The required condition for existence of decomposition $ \det 
\frac{\partial (\varphi_i, \ldots, \varphi_n)}{\partial (x_i, \ldots, x_n)}(a)
\neq 0$ for all $1\leq i \leq n$ is satisfied in an open dense subset of
$\mathrm{Dif}(U)$.
 
\eop

\noindent \textbf{Remark 3.} The statement of the lemma above has considered 
the maximal flag manifold in order to provide the total cascade decomposition
for flows in $\R^n$. The result can easily be restated in terms of just a pair
of complementary distribution of any dimension, as in the hipothesis of
the previous section. In the formulae for the components of the decomposition of
a flow $\varphi_t=\xi_t \circ \psi_t$ as in Theorem \ref{teorema decomposicao
global}, the equation of $\xi_t$ explodes when $\varphi_t$ reaches points
outside the dense open set containing the identity, as stated in the item (iv)
of Lemma \ref{decomposicao total}. An illustrative example is the linear planar
systems $(x', y')=(-y,x)$, where the rotation $\varphi_t$ hits a point
outside the dense open set. In fact, easy calculation of the components of
$\varphi_t=\xi_t \circ \psi_t$ shows that the entry $[\xi_t]_{1,2}$
satisfies 
\[
 \frac{d}{dt}[\xi_t]_{1,2} = -1 -[\xi_t]_{1,2}^2,
\]
 which explodes at $t=\pi/2$.

\end{document}